\documentclass [12pt]{amsart}

\usepackage{amssymb,amsxtra,amsfonts}
\usepackage{amsmath, accents}
\usepackage{graphicx}

\usepackage{graphics}

\usepackage{verbatim}
\usepackage{amssymb}
\usepackage{amsbsy}
\usepackage{amscd}
\usepackage{amsmath}
\usepackage{amsthm}
\usepackage[mathscr]{eucal}

\usepackage{verbatim}  
\usepackage{bbding}   

\newenvironment{ppb}[1]
{\ \!\!\!\!\!\!\!\!\!\!\!\!\!\!\!\!\!\!\!\!\!\!\!\!\!\!\!\!\!\!\!\!\!\!\!\!\!\!\!\! {\bf PPB------------------------------------------------------------------------------------------------PPB}\newline \tiny {#1}
\  \newline\normalsize\phantom{f}\!\!\!\!\!\!\!\!\!\!\!\!\!\!\!\!\!\!\!\!\!\!\!\!\!\!\!\!\!\!\!\!\!\!\!\!\!\!\!\! {\bf PPB------------------------------------------------------------------------------------------------PPB}\newline}{}

\long\def\pb #1*/{}

\def\reE@DeclareMathSymbol#1#2#3#4{%
    \let#1=\undefined
    \DeclareMathSymbol{#1}{#2}{#3}{#4}}
\DeclareSymbolFont{symbolsC}{U}{txsyc}{m}{n}
\SetSymbolFont{symbolsC}{bold}{U}{txsyc}{bx}{n}
\DeclareFontSubstitution{U}{txsyc}{m}{n}
\reE@DeclareMathSymbol{\strictiff}{\mathrel}{symbolsC}{76}

\newenvironment{ssn}[1]{{\ \newline \large\bf\hspace{-\parindent}{#1.}}}{}

\newcommand\beq{\begin{equation}}
\newcommand\eeq{\end{equation}}
\newcommand\bal{\begin{align*}}
\newcommand\eal{\end{align*}}   
\newcommand\bmx{\left(\begin{matrix}}
\newcommand\emx{\end{matrix}\right)}
\newcommand\bsmx{\left(\begin{smallmatrix}}
\newcommand\esmx{\end{smallmatrix}\right)}

\newcommand{\spq}{/\!\!/}

\providecommand{\spqa}[1]{\underset{#1}{/\!\!/}}

\newcommand{\st}{\ \bigl\vert\ }

\def\part#1{\frac{\partial\phantom{q}}{\partial#1}}

\newcommand {\flp}{(\!(}
\newcommand {\frp}{)\!)}

\newcommand{\fusion}[1]{\underset{#1}{\circledast}}

\newcommand{\DR}{\text{\rm DR}}
\newcommand{\Betti}{\text{\rm B}}
\newcommand{\MB}{\mathcal{M}_{\text{\rm B}}}

\newcommand{\Id}{\text{\rm Id}}

\newcommand{\Lie}{{\mathop{\rm Lie}}}

\DeclareMathOperator{\ISto}{{\IS}to} 

\DeclareMathOperator{\Sect}{\mathop{\rm Sect}}

\newcommand{\papk}[3]{\,_{#1}^{\phantom{#3}}\cA_{#2}^{#3}     }

\newcommand{\gatk}{\papk{G}{T}{k}} 

\newcommand{\Ad}{{\mathop{\rm Ad}}}

\newcommand{\Prod}{\prod}

\DeclareMathOperator{\Hom}{Hom}         

\newcommand{\SL}{{\mathop{\rm SL}}}
\newcommand{\PSL}{{\mathop{\rm PSL}}}
\newcommand{\GL}{{\mathop{\rm GL}}}

\DeclareMathOperator{\Rep}{\rm Rep}
\DeclareMathOperator{\End}{End}
\newcommand{\diag}{{\mathop{\rm diag}}}

\newcommand{\fl}{{\mathop{\footnotesize\rm {f}l} }}   

\newcommand{\hk}{{hyperk\"ahler }}   

\newcommand{\bD}{{\bf D}}

\newcommand{\bp}{{\bf p}}

\newcommand{\bs}{{\bf S}}
\newcommand{\bS}{{\bf S}}

\newcommand{\IA}{\mathbb{A}}

\newcommand{\IC}{\mathbb{C}}

\newcommand{\IH}{\mathbb{H}}

\newcommand{\IL}{\mathbb{L}}

\newcommand{\IP}{\mathbb{P}}

\newcommand{\IS}{\mathbb{S}}

\newcommand{\IZ}{\mathbb{Z}}

\newcommand{\A}{\mathcal{A}}
\newcommand{\cA}{\mathcal{A}}
\newcommand{\cB}{\mathcal{B}}

\newcommand{\cC}{\mathcal{C}}

\newcommand{\cF}{\mathcal{F}}
\newcommand{\cG}{\mathcal{G}}

\newcommand{\ch}{\eta}     

\newcommand{\cL}{\mathcal{L}}
\newcommand{\M}{\mathcal{M}}
\newcommand{\cM}{\mathcal{M}}

\newcommand{\cN}{\mathcal{N}}
\newcommand{\cO}{\mathcal{O}}

\newcommand{\cP}{\mathcal{P}}

\newcommand{\cW}{\mathcal{W}}

\newcommand{\gM}{       \mathfrak{M}     }

\newcommand{\g}{       \mathfrak{g}     }

\newcommand{\lh}{\mathfrak{h}}

\newcommand{\wt}{\widetilde}

\newcommand{\wh}{\widehat}

\newcommand{\be}{\beta}
\newcommand{\ga}{\gamma}

\newcommand{\Ga}{\Gamma}

\newcommand{\la}{\lambda}
\newcommand{\La}{\Lambda}

\newcommand{\Si}{\Sigma}
\renewcommand{\th}{\theta}

\renewcommand{\bar}{\overline}

\makeatletter
 \newlength{\typesize}
\makeatother

\newlength{\vvoff}
\newlength{\hhoff}

\def\mapright#1{\smash{
        \mathop{\longrightarrow}\limits^{#1}}}

\def\mapdown#1{\Big\downarrow
        \rlap{$\vcenter{\hbox{$\scriptstyle#1$}}$}}

\def\underset#1#2{\ \smash{\mathop{ #2 }\limits_{#1}}\ }

\newcommand{\pf}{\begin{bpf}}

\newcommand{\pfms}{\begin{bpfms}}
\newcommand{\epf}{\end{bpf}\hfill$\square$\\}           
\newcommand{\epfms}{\end{bpfms}\hfill$\square$\\}       

\newcommand{\idea}{\begin{bidea}}

\newcommand{\eidea}{\end{bidea}\hfill$\square$\\}           

\newcommand{\sk}{\begin{bsk}}    

\newcommand{\esk}{\end{bsk}\hfill$\square$\\}           
\newcommand{\sketch}{\begin{bsketch}}

\newcommand{\esketch}{\end{bsketch}\hfill$\square$\\}

\newtheorem {hypo}{\bf\hspace{-\parindent}Hypothesis}
\newtheorem {thm}[hypo]{Theorem}   
\newtheorem {prop}[hypo]{Proposition}

\newtheorem {cor}[hypo]{Corollary}
\newtheorem {lem}[hypo]{Lemma}

\newtheorem {defn}[hypo]{Definition}

\theoremstyle{remark}\newtheorem{rmk}[hypo]{Remark}

\renewenvironment{ppb}[1]{}{}

\begin{document}

\title[Wild Character Varieties, points on $\IP^1$ and Calabi's examples]{Wild Character Varieties, points on the Riemann sphere
and Calabi's examples}
\author{Philip Boalch}%

\maketitle

\section{Introduction}

Wild character varieties are moduli spaces of monodromy data of connections on bundles on smooth algebraic curves.
They were shown to admit holomorphic symplectic structures in 
\cite{thesis, smid}, to admit (complete) \hk metrics in \cite{wnabh}, and to arise as finite dimensional quasi-Hamiltonian quotients in \cite{saqh02, saqh, fission, gbs}.
A simple example was  shown to underlie the Drinfeld--Jimbo quantum group in \cite{smapg} (as conjectured in \cite{thesis, smid}) 
and further it was shown in \cite{bafi} that Lusztig's symmetries 
(a.k.a. the quantum Weyl group generators)
are the quantisation of  a simple  example of a 
wild mapping class group action on a wild character variety.

The wild character varieties  generalise the tame character varieties, which are moduli spaces of monodromy data of regular singular connections,
i.e. spaces of representations of the fundamental group. 
The extra monodromy data, enriching the fundamental group representation, needed to classify irregular connections is known as ``Stokes data''.
There are at least two approaches to Stokes data.
One approach ``Stokes filtrations''  (due to Deligne \cite{deligne78, bertrand-bbk79, Mal-irregENSlong, BV89,  malg-book, DMR-ci}, building on work of Malgrange, Sibuya and others)
involves adding flags on sectors at each pole, measuring the possible exponential 
growth rates of solutions.
In general it is complicated to classify such flags.
A different, but algebraically equivalent, approach was developed by Balser--Jurkat--Lutz
\cite{jurkat78, BJL79}, Martinet--Ramis \cite{MR91}, Loday-Richaud \cite{L-R94} and others. It involves canonical Stokes matrices and leads to the notion of ``Stokes local system''.  This approach was extended to arbitrary reductive groups $G$ in \cite{bafi}
and used in the description of the wild character varieties
as multiplicative symplectic quotients.

The aim of this article is to describe a simple class of examples of wild character varieties (studied in depth by Sibuya \cite{sibuya1975}) from both points of view, to illustrate this dichotomy.
In these examples it is not so difficult to directly bridge the gap between the two viewpoints.
A key point is that for connections on rank two bundles the flags amount to points of the Riemann sphere, and so in simple cases the wild character varieties are specific moduli spaces of points on $\IP^1$, studied by Sibuya when he considered the distinguished ``subdominant'' solutions in sectors at the poles.

For example we will explain the following theorem.
Let
$$\cM_{2k}^\text{Sibuya}  = \{ p_1,\ldots, p_{2k}\in \IP^1(\IC) \st 
p_1\neq p_{2}\neq \cdots \neq p_{2k}\neq p_1 \}/\PSL_2(\IC)$$
be the moduli space of $2k$-tuples of points of the Riemann sphere such that cyclically-consecutive points are distinct.
The prescription
$$\varphi(p_1,\ldots, p_{2k}) 
=  (-1)^k 
\frac
{(p_1-p_2)(p_3-p_4)\cdots (p_{2k-1} - p_{2k})}
{(p_2-p_3)(p_4-p_5)\cdots (p_{2k} - p_{1})}
$$
gives a well-defined map $\varphi:\cM_{2k}^\text{Sibuya}\to \IC^*$, generalising the cross-ratio, and the subvarieties $\cM_{2k}^\text{Sibuya}(q):=\varphi^{-1}(q)$ have dimension $2k-4$ for any $q\in\IC^*$.
\begin{thm}
For any $q\in \IC^*$ the space $\cM_{2k}^\text{Sibuya}(q)$ is a wild character variety
(by \cite{sibuya1975}), and  it is complex symplectic (by \cite{smid}). 
If $q\neq 1$ it is smooth, and it admits a complete \hk metric (by \cite{wnabh}). 
\end{thm}

If $k=3$ (i.e. $6$-tuples of points) then each space $\cM_{2k}^\text{Sibuya}(q)$ has real dimension four. Physicists refer to such complete \hk manifolds as ``gravitational instantons''. However these examples were not known to physicists. 
On the other hand the underlying complex algebraic surfaces have appeared frequently, in relation to the second Painlev\'e equation, as will be explained. 

The layout of this article is as follows.
Section 2 gives the direct ``canonical Stokes matrices'' approach to a simple class of  wild character varieties, and explains how they arise as finite dimensional multiplicative symplectic quotients. 
Section 3 then recalls the direct approach to the spaces of points on $\IP^1$ 
considered by Sibuya and then gives the direct proof of the 
isomorphism between the two  approaches.
(A more sophisticated approach is discussed in the appendix, which provides a brief introduction to Stokes filtrations and Stokes local systems.)
Next Section 4 describes the quiver approach and shows that these wild character varieties are multiplicative analogues of a family of \hk manifolds introduced by Calabi.
Finally Section 5 relates these examples to a 1764 paper of Euler, and
shows that Euler's continuant polynomials are group valued moment maps.
\section{Abelian Fission Spaces and Wild Character Varieties}

The quasi-Hamiltonian approach involves constructing the wild character varieties as finite dimensional {\em multiplicative} symplectic quotients. 
The symplectic/Poisson structure on the wild character variety 
then arises algebraically from the quasi-Hamiltonian two-form upstairs.
For the present examples only some simple quasi-Hamiltonian spaces will be needed.

Let  $G=\GL_2(\IC)$
and  consider the following subgroups:
$$ U_+ = \bmx 1 & * \\ 0 & 1 \emx \subset 
B_+ = \bmx * & * \\ 0 & * \emx 
\subset \ G\  \supset  
B_- = \bmx * & 0 \\ * & * \emx
\supset  
U_- = \bmx 1 & 0 \\ * & 1 \emx,
$$
and  
$T = B_-\cap B_+$ %
 the diagonal subgroup. 
The wild character varieties we will consider are as follows.
Choose an integer $k\ge 1$ and 
consider the variety 
\beq\label{eq: pwcv}
\MB = \{ 
\bS \in (U_+\times U_-)^k\st
 S_{2k}\cdots S_2S_1 \in T \} / T
 \eeq
 where 
$T$ acts by diagonal conjugation, and we take the affine geometric invariant theory quotient 
(the affine variety associated to the ring of $T$-invariant 
 functions).
 Here 
 $\bS =  (S_1,\ldots,S_{2k})$ with
$S_{\text{even}}\in U_-$ and $S_{\text{odd}}\in U_+$.
Further, for any fixed $t\in T$ of determinant one, consider the subvariety
\beq\label{eq: slsdescrn}
\MB(t) = \{ 
\bS \in (U_+\times U_-)^k\st
 S_{2k}\cdots S_2S_1 = t \} / T \subset \MB
 \eeq
which  has dimension $2k-4$.
It is a hypersurface in $\MB$.

To obtain these wild character varieties as multiplicative symplectic quotients consider the smooth affine variety
$$\cA = \gatk := G\times (U_+\times U_-)^k \times T.$$

In this section we will explain and illustrate the following result.

\begin{thm}[\cite{saqh02}\footnote{In fact \cite{saqh02} proves this for arbitrary complex reductive groups $G$ (with $B_\pm$ opposite Borels). The spaces $\cA$ are denoted $\wt \cC/L$ in \cite{saqh02} Rmk 4 p.6.
The non-abelian extension appears in \cite{fission, gbs}
(with $B_\pm$ replaced by arbitrary opposite parabolics, and $T$ by their common Levi subgroup).}]\label{thm: ab fission space}
The abelian fission space $\cA = \gatk$ is an algebraic quasi-Hamiltonian $G\times T$-space.
\end{thm}

This result means there is an action of $G\times T$ on $\cA$, an invariant algebraic two-form $\omega$ on $\cA$, and a {\em group-valued} moment map
$$\mu = (\mu_G,\mu_T): \cA \to G\times T;$$ 
$$\mu_G(C,\bs,h) = C^{-1} h S_{2k}\cdots S_2 S_1 C \in G,$$
$$\mu_T(C,\bs,h) = h^{-1} \in T$$
satisfying various axioms, which are multiplicative analogues of the usual axioms for a Hamiltonian $G\times T$-space.
The formula for $\omega$ is deferred to 
Remark \ref{rmk: twoform} below.
In particular  we will consider the reduction by $G$
$$\cB = \cB^k:= \gatk\spq G = \mu_G^{-1}(1)/G,$$
which is a smooth affine variety of dimension $2k-2$, 
with a residual  action of $T$.
Some immediate consequences of Theorem \ref{thm: ab fission space}, and the general quasi-Hamiltonian/quasi-Poisson yoga \cite{AMM, AKM},  are as follows.

\begin{cor} 1) The quotient 
$$\cA/G\cong (U_+\times U_-)^k\times T$$ 
is an algebraic Poisson manifold
 with symplectic leaves %
 $$\cA\spqa{\cC} G = \mu_G^{-1}(\cC)/G$$ 

\noindent
 for conjugacy classes 
 $\cC\subset G$,

 2) $\cB = \cA\spq G$ is an algebraic symplectic manifold,

3) The quotient by $T$ 
$$\mu_G^{-1}(\cC)/(G\times T)$$
of any symplectic leaf from 1) %
is a Poisson variety with symplectic leaves  
$\mu^{-1}(\cC\times\{t\})/(G\times T),$ for elements $t\in T$,

 4) For any $t\in T$ the reduction $\cB\spq_{\! t}\, T \cong \mu^{-1}(\{1\}\times\{t\})/G\times T$ 
 is a symplectic variety. 

\end{cor}
\pf
1) is a general result about quasi-Poisson manifolds; the quotient is a quasi-Poisson $T$-space, which means it is Poisson since $T$ is abelian. The leaves are the quasi-Hamiltonian reductions by $G$ (and are a priori quasi-Hamiltonian $T$-spaces, but this implies they are symplectic since $T$ is abelian). 
2) is a special case of 1), taking $\cC=\{1\}\subset G$.
3) follows the same pattern as in 1), considering the full action of $G\times T$. 
4) is a special case of 3).
\epf

Now it is easy to obtain the wild character varieties from the fission spaces.
The action of $G\times T$ on $\cA$ is given by 
$$(g,t)(C,\bS,h)  = (tCg^{-1}, t\bS t^{-1}, h)$$
where $(g,t)\in G\times T$ and 
$t\bS t^{-1} = (tS_1t^{-1},\ldots,tS_{2k}t^{-1})$.
Thus 
\begin{align}
\cB = \cA\spq G &= 
\mu_G^{-1}(1)/G \notag \\
&= \{ 
(\bS,h) \in (U_+\times U_-)^k\times T\st
 hS_{2k}\cdots S_2S_1 =1  \} \label{eq: bdescrn}  \\
&\cong \{ 
\bS \in (U_+\times U_-)^k\st
 S_{2k}\cdots S_2S_1 \in T \} \notag
 \end{align}
 which is a quasi-Hamiltonian $T$-space with moment map $h^{-1}$, and in particular a symplectic manifold (as $T$ is abelian).
 In turn  $\MB = \cB/T$ is thus a Poisson variety, with symplectic leaves 
 \beq\label{eq: mbt spq}
 \MB(t) = \cB\spqa{t} T = \{(\bS,h)\in \cB \st h^{-1}=t\}/T.
\eeq
Equivalently if $\cC\subset G\times T$ is the conjugacy class containing $1\times t\in G\times T$
then $$\MB(t) = \cA\spqa{\cC} G\times T = \mu^{-1}(\cC)/G\times T.$$

As mentioned above, the general quasi-Hamiltonian yoga 
implies that $\MB(t)$ is a symplectic variety: the restriction of $\omega$ to $\mu^{-1}(\cC)\subset \cA$ 
descends to give the symplectic form on $\MB(t)$.

The standard examples \cite{AMM} of quasi-Hamiltonian spaces are 
the conjugacy classes $\cC\subset G$ (with the inclusion being the moment map), and the double $\bD\cong G\times G$, which are multiplicative analogues of the coadjoint orbits $\cO\subset \g^*$ and the cotangent bundle $T^*G$ in the usual Hamiltonian world.
Using the fusion and reduction processes these examples give a clean algebraic construction of the Atiyah--Bott symplectic/Poisson structure on spaces of fundamental group representations of Riemann surfaces
\cite{AMM}.
The fission spaces give  an algebraic construction of the more general irregular Atiyah--Bott symplectic/Poisson structure of \cite{thesis, smid} on spaces of monodromy/Stokes data. 
The holomorphic symplectic manifolds in 3),4) above are examples of this general construction.
In these examples the symplectic form  was then computed explicitly 
by Woodhouse \cite{Woodh-sat} and this led to the 
above quasi-Hamiltonian approach (which works for arbitrary numbers of poles on arbitrary genus Riemann surfaces).
In turn the  irregular Atiyah--Bott complex symplectic quotients 
were upgraded into hyperk\"ahler quotients in \cite{wnabh}. 
In the present context this implies:

\begin{thm}[\cite{wnabh}]
If $t\neq 1$ then $\MB(t)$ is a complete hyperk\"ahler manifold of real dimension $4k-8$.  
\end{thm}

The condition $t\neq 1$ implies that all the points are stable, and so the spaces are smooth 
(cf. \cite{wnabh} \S8).
 In the set-up of \cite{smid, wnabh} this class of examples appears by considering 
 meromorphic connections on rank two bundles on the Riemann sphere, 
 having just one pole of order $k+1$.
 We then use the irregular Riemann--Hilbert correspondence to pass to the wild character variety.
(In general the \hk metric will depend on the choice of irregular type of the connection at the pole---cf. also 
\cite{ihptalk}). 
Note that in the context of compact K\"ahler manifolds, it is known that any holomorphic symplectic manifold admits a \hk metric, but the situation is more subtle in the noncompact case.
In the case $k=3$, $\MB(t)$ is a complete hyperk\"ahler manifold of real dimension four, i.e. a ``gravitational instanton'' in the physics terminology.
We will see below, when discussing quivers, that they may be viewed as multiplicative analogues of the Eguchi--Hanson spaces 
(the $A_1$ ALE spaces).

\begin{rmk} \label{rmk: FN cubic surface}
In the case $k=3$ it is easy to describe  the complex surface $\MB(t)$ explicitly.
It is isomorphic to the affine cubic surface 
\beq\label{eq: FNeq}
x\,y\,z + x + y + z  = b - b^{-1}
\eeq
where $b$ is a non-zero constant such that $t=\diag(q^{-1},q)$ and $q=-b^2$.
These cubic surfaces appear in \cite{FN80} (3.24).
In more detail let $s_i$ be the nontrivial off-diagonal matrix entry 
of $S_i$.
Then (cf.  \eqref{eq: slsdescrn}) the equation $S_6 \cdots S_1=t$ is equivalent to the three equations:
$s_1 = -q(s_3 s_4 s_5 + s_3 + s_5)$, 
$s_6 = -q(s_2 s_3 s_4 + s_2 + s_4)$
(allowing to eliminate $s_1$ and $s_6$), and 
 ${ s_2}{ s_3}{ s_4}{ s_5}+{ s_2}{ s_3}+{ s_2}{ s_5}+{ s_4}{ s_5}+1  = 1/q$.
 To quotient by $T$ we pass to invariants 
 $s_{23}, s_{25},  s_{34}, s_{45}$ where $s_{ij}=s_is_j$, and thus find
 $s_{25}=1/q-(1+s_{45}+s_{23}+s_{23}s_{45})$.
 Substituting this in the relation 
 $s_{23}s_{45}  = s_{34}s_{25}$
 yields the Flaschka--Newell surface \eqref{eq: FNeq} 
 after relabelling
$s_{45}=x/b-1, 
 s_{23}=y/b-1$ and $s_{34}=-1-b\,z$.
Note that Flaschka--Newell find these cubic surfaces as wild character varieties in a slightly different context (using a non-standard Lax pair for the Painlev\'e 2 equation), but their approach is known to be equivalent to the usual approach (cf. \cite{JKT09}), explaining why the cubic surfaces are the same. See \S\ref{sn: euler} for higher $k$. 
\end{rmk}

\begin{rmk} \label{rmk: twoform}
An explicit formula for the algebraic two-form $\omega$ on $\cA$ is as follows.
Define maps $C_i: \cA \to G$ by
$$C_i = S_i\cdots S_2S_1C$$
so that $C=C_0$, and let $b=hS_{2k}\cdots S_2S_1$.
Then
\beq  \label{eq: omeg def}
2 \omega = 
(\bar\ga,\Ad_b\bar\ga) 
+ (\bar\ga,\bar\be)
+ (\bar\ga_{2k}, \ch) -
\sum_{i=1}^{2k}(\ga_{i},\ga_{i-1})
\eeq
where the brackets $(\ ,\ )$ 
denote a fixed bilinear form on $\g$
and the Greek letters denote the following Lie algebra valued one-forms: 
$$\ga_i = C_i^*(\th),\qquad\bar\ga_i = C_i^*(\bar\th),\qquad
\eta = h^*(\th_T),\qquad
\bar\be = b^*(\bar\th)$$
where 
$\th=g^{-1}dg, \bar \th = dgg^{-1}$ are the left and right invariant Maurer--Cartan forms respectively. This expression equals that in \cite{saqh02}.
\end{rmk}

\section{Sibuya spaces}

Sibuya \cite{sibuya1975} studied the Stokes data at $\infty$ of 
differential equations of the form
\beq\label{eq: ode}
\frac{d^2y}{dz^2} = p(z) y
\eeq
for complex monic polynomials $p(z)$ of degree $m$.
In brief there are $m+2$ distinguished directions at $\infty$ (``Stokes directions'') and a preferred solution $v_i$ (the subdominant solution) on each sector between two consecutive Stokes directions. 
More precisely only the ray $\langle v_i \rangle$ spanned by $v_i$ is canonically determined.
Since the (rank two) local system of solutions on the plane is trivial, this determines 
$m+2$ rays in $\IC^2$, i.e. $m+2$ points of $\IP^1$. 
Sibuya thus considered the following spaces and related them to Stokes matrices. 

Let $n=m+2$ and let
$$X_{n}^\text{Sibuya} = \{ p_1,\ldots, p_{n}\in \IP^1(\IC) \st 
p_1\neq p_{2}\neq \cdots \neq p_{n}\neq p_1 \}$$
be the configuration space of $n$-tuples of cyclically ordered points of $\IP^1$ such that consecutive points are distinct.
Let
$$\cM_{n}^\text{Sibuya}  = X_{n}^\text{Sibuya}/\PSL_2(\IC)$$
be the geometric invariant theory quotient,
which has dimension $n-3$.
Note that there is an inclusion
$$\gM_{0,n} \ \subset \  \cM_{n}^\text{Sibuya}$$%
where $\gM_{0,n}$
is the moduli space of genus zero curves with $n$ distinct marked points.%

Suppose now that $n$ is even and set $n=2k$.
Consider the  function
$$\varphi : X_{2k}^\text{Sibuya} \to \IC^*$$
defined by
$$\varphi(p_1,\ldots, p_{2k}) 
=  (-1)^k 
\frac
{(p_1-p_2)(p_3-p_4)\cdots (p_{2k-1} - p_{2k})}
{(p_2-p_3)(p_4-p_5)\cdots (p_{2k} - p_{1})}.
$$
Up to a sign  this is the ``multiratio'' of the $2k$ points 
(see e.g. \cite{king-schief}). 
It is invariant under diagonal M\"obius transformations.

Choose an element $q\in \IC^*$ and let
$$\cM_{2k}^\text{Sibuya}(q) = \varphi^{-1}(q)/\PSL_2(\IC)\subset \cM_{2k}^\text{Sibuya}$$
be the subvariety of points with fixed multiratio.
It has dimension $2k-4$.

The stability condition for tuples of points on $\IP^1$ is well-known
(\cite{mumford.icm62}) and leads to the following

\begin{lem}
If $q\neq 1$ then all the points are stable, and so $\cM_{2k}^\text{Sibuya}(q)$ is the set of $\PSL_2(\IC)$-orbits in 
$\varphi^{-1}(q)\subset X_{2k}^\text{Sibuya}$.
\end{lem}
\pf
A $2k$ tuple $\bp$ is stable if no point has multiplicity $k$ or more. In the current set-up (with consecutive points distinct)  the multiplicity is at most $k$. Clearly if $k$ of the points are equal (e.g. $p_{\text{odd}} = 0$) then $\varphi(\bp)=1$.
\epf

\begin{rmk}
Note that if $k=3$ then the condition $\varphi(\bp)=1$ means that the multiratio of the $6$ points is $-1$, which is a condition  much-studied, even classically 
(see \cite{king-schief} and references therein).
\end{rmk}

Our aim now is to explain the following (which follows easily from \cite{sibuya1975}).
\begin{prop}\label{prop: sib.ss-sls}
The Sibuya moduli space $\cM_{2k}^\text{Sibuya}$ is algebraically isomorphic to the Poisson wild character variety $\MB$ in 
\eqref{eq: pwcv}, 
and the multi-ratio function cuts out the symplectic leaves:
$\cM_{2k}^\text{Sibuya}(q)$ is algebraically isomorphic to the wild character variety $\MB(t)$, where $t=\diag(q^{-1}, q)$.
\end{prop}

As discussed in the previous section, standard results on wild character varieties then imply that 
$\cM_{2k}^\text{Sibuya}(q)$ is a symplectic variety,  and even 
a complete \hk manifold whenever $q\neq 1$.

\ppb{
In fact the isomorphism $\cM_{2k}^\text{Sibuya} \cong \MB$ is a special case of a much more general result relating the two main ways to think about Stokes data (Stokes filtrations and Stokes local systems), so we will first explain the general context before returning to the proof of Proposition \ref{prop: sib.ss-sls}.
}

\pf
Write $n=2k$ and consider the affine variety 
$$V=V_{n}:=\{ v_1,\ldots, v_{n}\in \IC^2\setminus \{0\}
 \st 
v_1\nparallel v_{2}\nparallel \cdots \nparallel v_{n}\nparallel v_1 \}$$
where the symbol 
$\nparallel$ means ``is not parallel to''.
The torus $(\IC^*)^{n}$ acts on $V$ by scaling the vectors, and $G:=\GL_2(\IC)$ acts diagonally.
The quotient is $\cM_n^{\text{Sibuya}}$.
Viewing the $v_i$ as column vectors define  $2\times 2$ matrices 
$$ \Psi_i = (v_i v_{i+1})$$
(where the indices are taken modulo $n$).
By assumption these are all {\em invertible} matrices.
Thus following Sibuya we can define some matrices (``Stokes matrices''):
$$B_i:= \Psi_i^{-1}\Psi_{i-1}$$
so that  
$(v_{i-1} v_{i}) =  (v_{i} v_{i+1}) B_i$
(cf. \cite{sibuya1975} 21.32 p.86).
Clearly, by construction,
\beq\label{eq:sib prod}
B_{n}\cdots B_2B_1 = 1
\eeq
(cf. \cite{sibuya1975} 21.31) and moreover
$$B_i \in W:= \left\{\bmx * & 1 \\ * & 0 \emx\right\}\subset G $$
(as in \cite{sibuya1975} 21.30).
This leads to the space
$$\cW := \{B_1,\ldots,B_n\in W\st B_n\cdots B_2B_1=1\}$$
and the procedure above defines a map $\pi:V\to \cW$.

\begin{lem}\label{lem: trivbdle}
The map $\pi:V\to \cW$ is a trivial principal $G$-bundle.
In other words the action of $G$ on $V$ is free, has quotient $\cW$ and admits a global slice; $V\cong G\times \cW$. 
\end{lem}
\pf
Define an extended map $\wt \pi: V\to G\times \cW$ by setting $C=(v_1,v_2)=\Psi_1\in G$.
The formula $\Psi_i=\Psi_{i+1}B_{i+1}$ shows that each $\Psi_i$ is determined by $C=\Psi_1$ and the Stokes matrices, and in turn the $\Psi_i$ determine the $v_i$, so $\wt \pi$ is an isomorphism.
Specifically 
$\Psi_i = CB_1B_nB_{n-1}\cdots B_{i+1}$. 
The condition \eqref{eq:sib prod} ensures this holds for all $i$ modulo $n$.
\epf

This lemma holds even if $n$ is odd. 
To get to $\MB$ 
we adjust the 
ordering of some of the basis vectors to put the Stokes data in alternating Borels.
Suppose we swap the order of the columns of $\Psi_{2i}$, i.e. we redefine $\Psi_{2i}$ as 
$\Psi_{2i} = (v_{2i+1}, v_{2i})$ (and leave $\Psi_\text{odd}$ unchanged).
Then if we set $G_i = \Psi_i^{-1}\Psi_{i-1}$ we have 
$G_{2i}\in B_-, G_{2i+1}\in B_+$.
More precisely 
$$G_{2i}= PB_{2i} \in  W_-:=\left\{\bmx * & 0 \\ * & 1 \emx\right\},\quad
 $$$$G_{2i+1} = B_{2i+1}P \in
 W_+:=\left\{\bmx 1 & * \\ 0 & * \emx\right\}$$
 where
 $P=\bsmx 0 & 1 \\ 1 & 0 \esmx$.
This allows to rewrite Lemma \ref{lem: trivbdle}, showing that 
$$V/G =\cW \cong \cW':=\{G_1,\ldots,G_{2k}\in (W_+\times W_-)^k\st G_{2k}\cdots G_2G_1=1\}.$$

Now it is easy to compute the formal monodromy, and discover the multi-ratio.

\begin{lem}
The formal monodromy is $\diag(q,q^{-1})\in T$ where 
$q=\varphi(p_1,\ldots,p_{2k})$ is $(-1)^k$ times the multi-ratio of the points 
$p_i=\langle v_i\rangle \in\IP^1$.
\end{lem}
\pf
The formal monodromy is the monodromy of the associated graded local system (as in \cite{deligne78}), which here means we ignore the off-diagonal entries of the $G_i$ and only consider their diagonal  entries.
If we set $d(i,j) = \det(v_iv_j)$ then these nontrivial diagonal entries are
$$\det(G_{i}) = \frac{d(i-1,i)}{d(i+1,i)}$$
so that the formal monodromy is $\diag(q,q^{-1})$ where
\begin{alignat*}{1}
q&= \det(G_{2k}) \det(G_{2k-2})\cdots  \det(G_{2})
=\frac
{d(2k-1,2k)\cdots d(3,4)d(1,2)}
{d(1,2k)\cdots d(5,4)d(3,2)}\\
&=(-1)^k\frac
{(p_1-p_2)(p_3-p_4)\cdots (p_{2k-1}-p_{2k})}
{(p_2-p_3)(p_4-p_5)\cdots (p_{2k}-p_1)} =\varphi(p_1,\ldots,p_{2k}).
\end{alignat*}
\epf

Finally we need to consider the induced action of the torus $\wt T:=(\IC^*)^{2k}$ on $\cW'$ 
and show it reduces to an action of $T\cong (\IC^*)^2$ as in the definition of $\MB$.

Write an element of $\wt T$ as $c=(c_1,\ldots,c_{2k})$. This acts by scaling the $v_i$, and the induced action on the $G_i$ is:
$$c\cdot G_{2i} = \diag(c_{2i+1}, c_{2i}) \, G_{2i}\, \diag(c_{2i-1}, c_{2i})^{-1},$$ 
$$c\cdot G_{2i+1} = \diag(c_{2i+1}, c_{2i+2})\,  G_{2i+1}\, \diag(c_{2i+1}, c_{2i})^{-1}.$$ 
Thus in all cases the action on the nontrivial diagonal entry of $G_i$ is by multiplication by $c_{i+1}/c_{i-1}$.
We can thus use up most of this torus action by setting most of the $G_i$ to be unipotent:
Namely we can restrict to the subtorus $T\cong (\IC^*)^2\subset \wt T$ where 
$$c_1=c_3 =\cdots = c_{2k-1}\ \text{and}\ c_2=c_4 =\cdots = c_{2k}$$
and restrict to the subset $\cW''\subset \cW'$ where 
$G_2, G_3,\cdots,G_{2k-1}$ are unipotent (i.e. $\det(G_i)=1$ if $i\neq 1,2k$).
Then the quotient $\cW'/\wt T$ is identified with the quotient $\cW''/T$.
If we write 
 $$G_1 = S_1\diag(1,q^{-1})$$
 with $S_1\in U_+$, %
 and
 $$G_{2k} =  \diag(q,1)S_{2k}$$
 with $S_{2k}\in U_-$
and set $S_i=G_i\in U_\pm$ otherwise, then the relation $G_{2k}\cdots G_1=1$
becomes
$$h S_{2k}\cdots S_2S_1 = 1$$
where $h=\diag(q,q^{-1})$ is the formal monodromy, so that 
$$\cW''\cong \{ \bS\in (U_+\times U_-)^k\st S_{2k}\cdots S_1\in T\} = \cB.$$
Moreover the $T$-action on $\cW''$ matches the $T$-action in the definition \eqref{eq: pwcv} 
of $\MB$, and so the proof is complete.
\epf

\begin{rmk}
This algebraic isomorphism is an example of the equivalence of categories between Stokes filtered local systems and Stokes local systems (both of which are equivalent to a category of connections, cf. \cite{deligne78}, \cite{gbs} Thm A.3).
It is a consequence of the simple analytic fact that
the columns of the canonical fundamental solutions (used to define the canonical Stokes matrices, as in \cite{BJL79} Thm A) are consecutive subdominant solutions (when suitably scaled and ordered), 
cf. Appendix A.
Loday-Richaud's algorithm \cite{L-R94} gives an algebraic procedure for translating between Stokes filtrations and Stokes local systems in general.
\end{rmk}

\begin{rmk}
Note that Sibuya related the spaces $\cM^\text{Sibuya}_r$ to Nevanlinna's theory of Riemann surfaces (see \cite{sibuya1975} p.ix).
\end{rmk}

\section{Multiplicative quiver varieties}

There is a standard way to associate \hk manifolds to graphs 
\cite{Kron.ale, nakaj-duke94}, known as additive/Nakajima quiver varieties.
Here we will briefly recall the algebraic 
approach to the underlying  
holomorphic symplectic varieties (cf. also \cite{cas-slod, CB-mmap})
and then discuss the multiplicative version relevant to the present examples.
Let $\Ga$ be a graph with nodes $I$.
Suppose we are given a vector space $V_i$ for each node $i\in I$.
Then we can consider the vector space 
$$\Rep(\Ga,V)=\bigoplus_{a\in \bar\Ga} \Hom(V_{t(a)},V_{h(a)})$$  
of maps along each edge of the graph in both directions.
We will call this the space of representations of the graph on the $I$-graded vector space 
$V=\bigoplus V_i$.
Here $\bar\Ga$ is the set of oriented edges of $\Ga$, i.e. the set of pairs $(e,o)$ where $e$ is an edge of $\Ga$ and $o$ is one of the two possible orientations of $e$.
Given an oriented edge $a\in \bar\Ga$, the head $h(a)\in I$ and 
tail $t(a)\in I$ are well-defined.

 The group $H:=\Prod\GL(V_i)$ acts on $\Rep(\Ga,V)$ via its natural action on $V$ preserving the grading.
Further a choice of orientation of the graph $\Ga$ determines a holomorphic symplectic structure on $\Rep(\Ga,V)$, and then the action of $H$ is Hamiltonian with a moment map
$$\breve\mu : \Rep(\Ga,V) \ \to\ \lh^*=\Lie(H)^*\cong \Prod_{i\in I}\End(V_i).$$

The %
Nakajima quiver varieties
 are defined by choosing a central value $\la\in \IC^I$ of the moment map and taking the symplectic quotient:
$$\cN(\Ga,d,\la) = \Rep(\Ga,V) \spqa{\la} H 
= \{ \rho\in \Rep(\Ga,V) \st \breve\mu(\rho) = \la\}/H$$
where  $\la$ is identified with the central element
$\sum\la_i\Id_{V_i}$ of $\Lie(H)^*$.
Here $d\in \IZ^I$ denotes the dimension vector, with components 
$d_i=\dim(V_i)$, and the spaces are empty unless $\sum \la_id_i=0$.
The quotient is the affine quotient, taking the variety associated to the ring of $H$ invariant functions (although often one adds an extra parameter choosing  a nontrivial linearisation as well---for simplicity here we won't do this).

\begin{figure}[ht]
	\centering
	\input{p2ngraph.pstex_t}
	\caption{The graph $\Ga$.}\label{fig: p2n graph}
\end{figure}

If $\Ga$ is an affine Dynkin graph and $d$ is the minimal 
imaginary root, then as discovered in \cite{Kron.ale},
$\cN(\Ga,d,\la)$ has complex dimension two,
and is a deformation of the Kleinian singularity 
$\IC^2/G(\Ga)$, where $G(\Ga)\subset \SL_2(\IC)$ 
is the McKay group of $\Ga$.
For type $A$, these \hk four-manifolds were known before: for $A_1$ they are the Eguchi-Hanson spaces \cite{EH78}
which in one complex structure are generic coadjoint orbits of 
$\SL_2(\IC)$ (and in another they are $T^*\IP^1$). 
Just after this example was found Calabi \cite{calabi-hk} found examples in all dimensions:  
In one complex structure Calabi's examples are minimal semisimple coadjoint orbits of 
$\SL_n(\IC)$ (and in another they are $T^*\IP^{n-1}$).
As quiver varieties, Calabi's examples arise by considering the 
simple graph in Figure \ref{fig: p2n graph}, with two nodes and $n$ edges (with dimension vector $d=(1,1)$).

In this example $\Rep(\Ga,V)$ has dimension $2n$ and $\cN$
is the symplectic reduction by a torus $H=(\IC^*)^2$, and has 
complex dimension $2n-2$ (the diagonal subgroup of $H$ acts trivially).
The aim of the rest of this section is to explain the following statement:

{\em The wild character varieties $\MB(t)$ in 
\eqref{eq: slsdescrn} are multiplicative analogues of Calabi's examples.}

First of all note that if we set $k=n+1$ then 
$\dim(\MB(t))= \dim(\cN(\Ga,d,\la))= 2n-2 = 2k-4$.

Secondly note that $\dim(\cB)= \dim(\Rep(\Ga,V)) = 2n$, and moreover:

\begin{lem} [cf. \cite{cmqv} Rmk 5.4]
The space $\cB=\cB^k$ may be identified with a Zariski open affine subset 
of  $\Rep(\Ga,V)$.
\end{lem}
\pf
In brief $V=V_1\oplus V_2$ with both $V_1,V_2$ one dimensional complex vector spaces,
and $\rho\in \Rep(\Ga,V)$ consists of $n$ maps $V_1\to V_2$ 
and $n$ maps $V_2\to V_1$.
Suppose we label these maps
$$s_2, s_4,\ldots s_{2n}:V_1\to V_2,\qquad s_1,s_3,\ldots s_{2n-1}:V_2\to V_1.$$
Then we may identify $\Rep(\Ga,V)$ with 
$(U_-\times U_+)^n\subset \GL(V)^{2n}$ by setting
$$
S_{2i} = \bmx 1 & 0 \\ s_{2i} & 1 \emx \in U_-\quad\text{and}\quad 
S_{2i+1} = \bmx 1 & s_{2i+1} \\ 0 & 1 \emx \in U_+.
$$
Now recall from \eqref{eq: bdescrn} that, with $k=n+1$
\begin{align}
\cB &= \{ 
\bS \in (U_+\times U_-)^k \st
 (S_{2n+2}S_{2n+1})S_{2n}\cdots S_2S_1 \in T  \}  \notag \\
&\cong \{ 
\bS \in (U_+\times U_-)^n\st
 S_{2n}\cdots S_2S_1 \in G^\circ \} \label{eq: bdescrn2}
 \end{align}
 where $G^\circ = U_+TU_-=U_+U_-T\subset G$ is the opposite Gauss cell.
 Thus \eqref{eq: bdescrn2} and the Gauss decomposition says that 
 $$\cB \cong \{\rho=(s_1,\ldots s_{2n}) \in \Rep(\Ga,V)\st
 (S_{2n}\cdots S_2S_1)_{22} \neq 0 \}$$
so that $\cB$ is indeed a Zariski open affine subset of $\Rep(\Ga,V).$
\epf

Now define the {\em invertible} representations of the graph $\Ga$ to be this open subset:
$$\Rep^*(\Ga,V) =  \{\rho  \in \Rep(\Ga,V)\st
 S_{2n}\cdots S_2S_1 \in G^\circ \} \subset \Rep(\Ga,V).$$
Thus, noting that $H=T$ is the maximal torus of $G=\GL(V)$, and that $\cB$ is a quasi-Hamiltonian $T$-space (cf. \eqref{eq: bdescrn}), there is a group-valued moment map
\beq \label{eq: qmm}
\mu:\Rep^*(\Ga,V)  \to H\eeq
defined by taking the $T$ component
of $S_{2n}\cdots S_2S_1 \in G^\circ$.
Beware $\mu$ is
 {\em not} the restriction of the usual moment map
$\breve\mu:\Rep(\Ga,V)  \to \Lie(H).$

In turn it is natural to define the multiplicative quiver varieties of $\Ga$ to be the multiplicative symplectic reductions of $\Rep^*(\Ga,V)$
at central values of the moment map.
Namely we choose $q\in (\IC^*)^I$ and define
$$\cM(\Ga,d,q) = \Rep^*(\Ga,V)\spqa{q} H = \mu^{-1}(q)/H.$$

\noindent
These will be empty unless $q^d:= \Prod q_i^{d_i}=1$, which in the current setup means $q_1q_2=1$.
Of course, this is just a rephrasing of the 
construction of the wild character variety $\MB(t)$ in 
\eqref{eq: mbt spq}, with $t=\diag(q_1,q_2)$, so that
$$\cM(\Ga,d,q) \cong \MB(t)$$
but now we see the link to graphs, and thus to Kac--Moody root systems.
(Many other examples appear in \cite{cmqv}.)
In the simplest nontrivial example $k=3$ the graph $\Ga$ is the affine $A_1$ Dynkin graph (with two edges)
and so the corresponding wild character varieties are multiplicative analogues of the Eguchi--Hanson spaces. 
Note that Okamoto 
\cite{Okamoto-dynkin} already related the second Painlev\'e equation to the affine $A_1$ Weyl group
(cf. also  \cite{quad} Exercise 3 and  \cite{rsode} Appendix C).

If we are prepared to work analytically then a more direct link between the additive and multiplicative quiver varieties is available.
Namely the Riemann--Hilbert--Birkhoff map plays a role analogous to the exponential map:

\begin{thm}
The Riemann--Hilbert--Birkhoff map, taking a connection to its monodromy data, gives a holomorphic map from the 
additive quiver varieties for the graph $\Ga$ to the corresponding multiplicative quiver varieties.
It relates the holomorphic symplectic structures (but not the \hk metrics).
\end{thm}
\pf
This follows from \cite{smid} Thm 6.1 modulo relating both sides  to quivers:
on the  multiplicative side this follows from the discussion above plus the quasi-Hamiltonian approach in \cite{saqh02} (see especially Cor. p.3).
The general dictionary relating the additive side to Nakajima quiver 
varieties was written down (and established in some cases) in \cite{rsode} Appendix C
and justified in general in \cite{hi-ya-nslcase}.

In more detail recall the  diagram of moduli spaces (from \cite{smid} (29) p.181):
\begin{equation}	\label{main diagram}
\begin{array}{cccccc}
&&  \wt\M_\DR & \mapright{\cong} & \wt\A_{\fl}/\cG_1 &\\
&&  \bigcup && \mapdown{\cong}&\\
\wt O_1\times\cdots\times\wt O_m\spq G &\cong&
  \wt\M_\DR^* & \mapright{\wt\nu} & \wt \M_\Betti. &\hphantom{\longrightarrow}
\end{array}
\end{equation}
Here $\wt\M_\DR^*$ is a moduli space of framed meromorphic connections on the trivial bundle on the Riemann sphere with $m$ poles, and $\wt\M_\Betti$ is the corresponding space of monodromy/Stokes data.
(See \cite{smid} for the full definitions.)
Thm 6.1 of \cite{smid} says that the map $\wt \nu$ taking monodromy/Stokes data is symplectic.
In turn in \cite{saqh02} the space $\wt\M_\Betti$ is 
identified (as a symplectic manifold) with a fusion product 
$$\wt\M_\Betti\cong \wt \cC_1\fusion{}\cdots \fusion{} \wt \cC_m\spq G$$
for certain quasi-Hamiltonian $G\times T$ spaces $\wt \cC_i$.
Thus the (framed) Riemann--Hilbert--Birkhoff map is a symplectic map 
$$\wt O_1\times\cdots\times\wt O_m\spq G \ \mapright{\wt\nu} \ 
\wt \cC_1\fusion{}\cdots \fusion{} \wt \cC_m\spq G.$$
It is injective and intertwines the $T$-actions changing the framings.
(This injectivity is due to \cite{JMU81}---it was extended to a bijective correspondence 
$\wt\M_\DR \cong \wt \cM_\Betti$ in \cite{smid} Cor. 4.9.) 
Now if we specialise to the case $m=1$ with just one pole and 
use \cite{smid} Lem. 2.4 to ``decouple''  
$$\wt O_1 \cong O_B\times T^*G$$
the space $\wt O_1$, then we obtain a $T$-equivariant symplectic map 
$$O_B \ \mapright{\wt\nu} \ \wt \cC_1\spq G.$$
Here $O_B$ is a coadjoint orbit of the (unipotent) group of based jets of bundle automorphisms, and thus has global Darboux coordinates (as noted in \cite{smid} p.190).
The conjecture of \cite{rsode} Appendix C proved in \cite{hi-ya-nslcase} 
identifies $O_B$ (as a Hamiltonian $T$-space) 
with a space of graph representations, and thus the reduction by $T$ with a Nakajima quiver variety.
If the group $G$ is $\GL_2$ then the graph is the graph $\Ga$ of Figure 
\ref{fig: p2n graph} if the connections have poles of order $k+1$.
On the other side we now just need to identify the $T$ reductions of 
$\wt \cC_1\spq G$ with the $T$ reductions of $\cB$, which is immediate from the definitions.
Indeed, the space $\wt \cC_1$ is a cover of the fission space $\cA$ used to define $\cB$.
\epf

\section{Eulerian hypersurfaces}\label{sn: euler}

On p.55 of Euler's 1764 article ``Specimen Algorithmi Singularis'' \cite{euler281}
the reader will find the list of polynomials:
\begin{align*}
(a) &= a\\
(a,b) &=  ab+1 \\
(a,b,c) &= abc+c+a \\
(a,b,c,d) &= abcd + cd + ad + ab +1 \\
(a,b,c,d,e) &= abcde + cde + ade + abe + abc + e +c +a \\
&\qquad\text{etc.}
\end{align*}
They appear when computing continued fractions and 
nowadays they are known as ``Euler's continuants''.
The latter monomials in the $n$th continuant arise by forgetting all possible pairs of consecutive letters from the first monomial.
Recalling Remark \ref{rmk: FN cubic surface} we thus see that, for $k=3$, the space
$\MB(t)$ is the quotient of the ``Eulerian hypersurface'' 
$$(s_2,s_3,s_4,s_5) = 1/q$$
by the action of $\IC^*$, and the expressions for $s_1,s_6$ involve continuants of degree  one less.
Moreover the quasi-Hamiltonian $T$-space $\cB$ is isomorphic to the open subset 
$$(s_2,s_3,s_4,s_5)\neq 0$$ of $\IC^4$ 
and the continuant $(s_2,s_3,s_4,s_5)$ is a group-valued moment map for the $\IC^*$ action on this subset.
The aim of this section is to point out that this all holds for any $k$:
\begin{prop}
For any $k$ the quasi-Hamiltonian $T$-space $\cB=\cB^k$ is isomorphic to the open subset 
$$(s_2,s_3,\ldots,s_{2k-1})\neq 0$$ of 
$\IC^{2k-2}$, 
the continuant $(s_2,s_3,\cdots,s_{2k-1})$ is a group-valued moment map for the $\IC^*$ action on this subset
and 
the wild character variety $\MB(t)$ is the quotient of the Eulerian hypersurface 
$$(s_2,s_3,\cdots,s_{2k-1}) = 1/q$$
by the action of $\IC^*$.
\end{prop}
\pf
The continuants can be defined (see e.g. \cite{knuth}) by the recurrence
$$(x_1,x_2,\ldots,x_n) = x_1(x_2,x_3,\ldots,x_n) + (x_3,x_4,\ldots,x_n)$$
and using this one can easily show 
\beq \label{eq: euler matrix}
\bmx x_1 & 1 \\ 1 & 0 \emx 
\bmx x_2 & 1 \\ 1 & 0 \emx \cdots
\bmx x_n & 1 \\ 1 & 0 \emx 
 = 
 \bmx 
 (x_1,\ldots,x_n) &  (x_1,\ldots, x_{n-1}) \\ 
 (x_2,\ldots,x_n) & (x_2,\ldots,x_{n-1}) \emx. 
\eeq
Now using the Gauss decomposition $\cB$ is isomorphic to 
$$\{  (S_2,\ldots,S_{2k-1}) \in (U_-\times U_+)^{k-1}\st
 (S_{2k-1}\cdots S_2)_{11} \neq 0 \}.$$
Thus the key point is to show that
$(S_{2k-1}\cdots S_2)_{11}  = (s_2,s_3,\ldots,s_{2k-1})$
where $s_i$ is the active off-diagonal entry of $S_i$.
But this is a direct consequence of \eqref{eq: euler matrix}, upon
noting that 
$P S_2 = \bsmx s_2 & 1 \\ 1 & 0 \esmx$  
and 
$ S_3 P  = \bsmx s_3 & 1 \\ 1 & 0 \esmx$  
etc., where $P=\bsmx 0 & 1 \\ 1 & 0 \esmx$.
\epf

\begin{rmk}
Although it is beyond the scope of this article let us mention that this gives a 
presentation of the  fission algebra $\cF^q(\Ga)$ of the graph $\Ga$: %
it is isomorphic to the quotient of the path algebra $\cP(\bar \Ga)$ of the quiver 
$\bar \Ga$ by the relations 
$$(a_1,b_1,a_2,b_2,\ldots,a_n,b_n)e_1 = q_1e_1$$
$$(b_n,a_n,\ldots,b_2,a_2,b_1,a_1)e_2 = q^{-1}_2e_2$$
where $q_i\in\IC^*$, the $a$'s are the arrows to the left, the $b$'s are the arrows to the right and 
$e_1/e_2$  is the idempotent for the left/right node respectively.
Here the ordering of the symbols in the continuants is as written by Euler. 
For $n=1$ one gets the multiplicative preprojective algebra $\La^q$ of \cite{CB-Shaw}
(which in this case is isomorphic to a deformed preprojective algebra), 
but in general $\cF^q\not\cong \La^q$ 
(cf. \cite{cmqv} Rmk 6.10).
\end{rmk}

\noindent{\bf Conclusion}

We have described four ways to think about a simple class of wild character varieties. 
In these examples the structure group was $\GL_2(\IC)$ and the underlying curve was the complex plane---the connections just had one singularity, at $\infty$,
and we assumed there was an even number of Stokes directions at the singularity.
Clearly these spaces admit lots of generalisations, and we hope they serve as a helpful introduction to the more general moduli spaces constructed in \cite{smid, saqh02, wnabh, fission, gbs}, complementing the survey \cite{p12}.
As mentioned above, Loday-Richaud's algorithm \cite{L-R94} gives an algebraic procedure for translating between Stokes filtrations and Stokes local systems in general.
Note that under this dictionary the spaces constructed in \cite{smid, saqh02} involve full flags, whereas the spaces constructed in \cite{wnabh, fission, gbs} involve arbitrary partial flags.
In a subsequent article we will describe in detail the ``odd'' case (also considered by Sibuya)---this is in some sense simpler since the spaces $\cM^\text{Sibuya}_{2k+1}$ are smooth symplectic varieties directly, with the same expression for the two-form $\omega$.
Moreover the analysis of \cite{wnabh} extends directly (as pointed out in \cite{sabbahrmk, witten-wild}) to show they are complete \hk manifolds.
In another direction a nice class of examples generalising the present ones come from the following graphs (bearing in mind the dictionary in \cite{rsode} Apx C, \cite{cmqv} \S3.3). 
In this case the additive quiver varieties are arbitrary coadjoint orbits of $\GL_n(\IC)$ (and in special complex structures they are cotangent bundles of flag varieties via \cite{nakaj-duke94} \S7, deframed in the sense of \cite{CB-mmap} p.261).

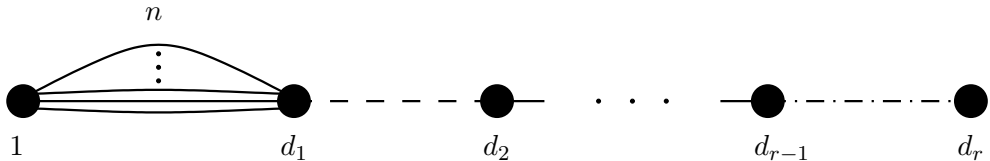
\begin{figure}[ht]
	\centering
	\input{Ographv2.pstex_t}
\end{figure}

\ 

\noindent{\bf Acknowledgments.}
Thanks are due to Alistair Scott MacLeod for providing a copy of  
\cite{knuth} at an opportune moment.

\ppb{
higher rank, coads.

(and we assumed it was untwisted there).

Thus we have described three ways to think about a simple class of wild character varieties. 
In this example the structure group was $\GL_2(\IC)$ and the underlying curve was the complex plane---the connections just had one singularity, at $\infty$ (and we assumed it was untwisted there).
Clearly this example admits lots of generalisations.
For example we can replace the complex plane with an arbitrary smooth algebraic curve and replace $\GL_2(\IC)$ by an arbitrary complex reductive group. If the singularities are generic then such moduli spaces are symplectic \cite{smid} and were shown to appear as multiplicative symplectic quotients in \cite{saqh02}.
In terms of flags this means we restrict to full flags and that consecutive flags are transverse (as you go around the circle of direction at each pole).
More generally one can have various configurations of partial flags and they become difficult to classify from the flag viewpoint. 
Indeed Malgrange \cite{Mal79} and Sibuya \cite{sibuya77} 
(see also \cite{BV89}) wrote down a non-abelian cohomology space which classifies the Stokes filtrations 
 \cite{deligne78}.
In turn  one can find canonical cocycles to describe this 
cohomology space \cite{L-R94}---but these cocycles  {\em are} (essentially)  the canonical Stokes multipliers.
Thus we could skip this discussion if we take the canonical Stokes multipliers as the basic objects and construct the wild character varieties in terms of them.
This leads to the notion of Stokes local system and is the approach taken in \cite{gbs}.
}

\appendix

\section{Stokes filtrations and Stokes local systems}

This appendix briefly recalls the notions of Stokes filtrations and Stokes local system in the present context.
This is more sophisticated than the presentation in the body of the text, and we will explain how the data there appears from the data here.
Suppose $(V,\nabla)$ is 
an algebraic connection on a rank two vector bundle $V$ on the complex plane, with unramified formal normal form at infinity.
In a global trivialisation such a connection takes the form
$$ \nabla = d-A$$
where $A$ is a $2\times 2$ matrix  of algebraic one-forms on $\IC$, i.e.
if $z$ is a coordinate on $\IC$ then  
$A = Bdz$
with $$B = \sum_0^m B_iz^i$$
a polynomial matrix.
Thus horizontal sections of $\nabla$ correspond to solutions of the differential 
equation $$\frac{dv}{dz} = B v$$
for vectors $v(z)$ of holomorphic functions 
on $\IC$.
For example a polynomial differential equation 
$y'' = p(z) y$ as in \eqref{eq: ode} 
determines the connection
with $B = \bsmx 0 & 1 \\ p & 0 \esmx$.
Changing the trivialisation of $V$ corresponds to the gauge action 
$A\mapsto g[A] = gAg^{-1} + (dg)g^{-1}$
for $g\in G(\IC[z])$, where $G=\GL_2$.
Restricting to the formal neighbourhood of $\infty$
amounts to considering the orbit of $A$ under the larger group $\wh G = G(\IC\flp z^{-1}\frp)$.
The condition of having unramified formal normal form at infinity means that 
the $\wh G$ orbit of $A$ contains an element of the form
$$A^0 = dQ + \La\frac{dz}{z}$$
for some element
$Q = \sum_1^k Q_iz^i$ where each $Q_i$ is a diagonal matrix.
We will suppose $Q$ is not central (i.e. at least one of the $Q_i$ is not a scalar matrix).
Then there is no loss of generality assuming $\La$ is a (constant) diagonal matrix. 
For example this occurs if the leading coefficient of $A$ is regular semisimple.
(In general one can take a root of $z$ to reduce any connection to one with unramifed formal normal form.)
The element $Q$ is the 
``irregular  type'' of $\nabla$.
Thus by assumption
$$A = \wh F [A^0]$$
for some $\wh F\in \wh G$.
Without loss of generality (tensoring by a connection on a line bundle) we may assume $Q_k$ is not a scalar matrix.
(Note that the leading coefficient of $A$ does not need to be semisimple---for example the polynomial differential equation \eqref{eq: ode} has unramified normal form if and only if $p$ has even degree.)
The connection $\nabla^0 = d-A^0$ 
has fundamental solution $z^\La e^Q$ and so 
$\nabla$ has {\em formal} fundamental solution 
$\wh F z^\La e^Q$.
Note that in general $\wh F$ will not be convergent.

\ssn{Stokes filtrations} 
This rests on the local asymptotic existence theorem, which
says that any real direction $d$ heading towards $\infty$ has a sectorial neighbourhood
$\Sect_d$ 
on which there is an analytic map
$F:\Sect_d\to G$ asymptotic to $\wh F$ with $F[A^0]=A$.
In general $F$ is not uniquely determined by these conditions. Nonetheless, given such $F$ we get a fundamental solution $\Phi_d=Fz^\La e^Q$ of $\nabla$ on $\Sect_d$ (which we can then freely analytically continue as a solution).
Now suppose we write $Q=\diag(q_1,q_2)$ with 
$q_i\in z\IC[z]$, $\La = \diag(\la_1,\la_2)$, and $F=(f_1,f_2)$ with $f_i:\Sect_d\to \IC^2$.
Let $v_i=f_iz^{\la_i}e^{q_i}$ be 
the columns of $\Phi_d$.
Then on $\Sect_d$ any horizontal section of 
$\nabla$ is a linear
combination of the horizontal 
sections $v_1$ and $v_2$.

Now the growth/decay of $v_i$ is dominated by the exponential factor $e^{q_i}$ so we should examine how this varies as a function of the direction $d$. Thinking about $e^z$ and then $\exp({z^k})$ shows there are $2k$  sectors where they alternately grow and decay.
This is nicely encoded in the {\em Stokes diagram} of the irregular type: draw a small dashed circle around 
$\infty$ and for each $i$ draw a wiggly curve around $\infty$ that is near $\infty$ when $\exp(q_i)$ grows and crosses the circle when 
$\exp(q_i)$ changes to exponential decay.
 In other words for each exponential factor and for each direction plot the possible growth rates, 
so the ordering of the intersections of these curves with any given direction corresponds to the dominance ordering.
A typical picture is as on the left of Figure \ref{fig: hexagon} (for the case $k=3$).
In the current setting the Stokes diagram has two components.
The Stokes diagram for the Airy equation appears in Stokes' paper \cite{stokes1857} p.116.
Define the {\em Stokes directions} to be the directions where the dominance ordering of the functions $e^{q_i}$ changes (i.e. where the wiggly curves cross).

\begin{figure}[ht]
	\centering
	\input{hexagon.pstex_t}
	\caption{}\label{fig: hexagon}
\end{figure}

Now suppose we choose a direction $d$ which is not a Stokes direction. Then the exponential factors have a definite dominance ordering along $d$
and  so one of the $v_i$ will be the least dominant solution (i.e. ``subdominant''). 
Suppose it is $v_1$ (for this $d$). 
Any other solution will be of the form $av_1 +b v_2$ so that  if $b\neq 0$ it will not be subdominant.
Thus the ray 
$\langle v_1 \rangle$ is uniquely determined by $\nabla$.  If we move $d$ slightly (so we don't cross a Stokes direction) then the same solution will remain subdominant.

The idea of Stokes filtered local systems is to axiomatise the rays that appear in this way  (and prove that the resulting topological objects classify the original connections):
Suppose we take $d$ to be a Stokes direction in the local asymptotic existence theorem.
Then $v_1$ will be the subdominant solution on one side of $d$ and $v_2$ will be the subdominant solution on the other side.
So the corresponding rays are tranverse (when continued across $d$).

Let $\Si=(\IP^1,\infty,Q)$ be the wild Riemann surface consisting of the Riemann sphere with the point $\infty$ marked and the irregular type $Q$ fixed.
Let $\Si^\circ = \IP^1\setminus \infty = \IC$ and draw the Stokes diagram of $Q$ on $\Si^\circ$.

\begin{defn}
A ``Stokes filtered local system'' 
for $\Si$ consists of:
a rank two local system $\cL\to \Si^\circ$ plus,
for each sector $\Sect_i\subset \Si^\circ$ at $\infty$ bounded by consecutive Stokes directions, a rank one sublocal system 
$\cL_i \subset \cL\bigl\vert_{\Sect_i}$
such that $\cL_i$ and $\cL_{i+1}$ are transverse 
(when extended across the Stokes direction between them).  
\end{defn}

More generally one takes the filtration of sublocal systems corresponding to the exponential growth rates of solutions (see \cite{deligne78, bertrand-bbk79, Mal-irregENSlong, BV89,  malg-book, DMR-ci}).
The associated graded local system corresponds to the formal normal form, and the tranversality condition may be expressed as saying that locally (around the circle of directions) the Stokes filtered local system is isomorphic to that determined by the associated graded (using the dominance ordering of the exponential factors).

Now it is easy to classify the Stokes filtered local systems that appear in the current setting. Note that $\cL$ is trivial (as $\IC$ is simply connected). 
Thus $H^0(\cL)$ is a two-dimensional complex vector space
and each $\cL_i$ extends uniquely to a global rank one subsystem of $\cL$ and so gives a one dimensional subspace
$H^0(\cL_i)\subset H^0(\cL)$.
Thus we get a $2$d vector space with a cyclically ordered collection of $1$d subspaces, such that consecutive subspaces are transverse, as studied in the body of this article following Sibuya.

As explained in \cite{BV89} in general it is not so easy to classify Stokes filtrations: in general one passes via the Malgrange--Sibuya nonabelian cohomology space and must parameterise that.
However it turns out there are preferred cocycles
and so a nice parameterisation is possible 
\cite{BJL79, MR91, L-R94}.
This leads to the notion of Stokes local system, which are  topological objects that are simpler to classify.
They may be approached directly (without first passing through Stokes filtrations) as follows.

\ssn{Stokes local systems}
This is an alternative way to extract topological data from an irregular connection and rests (in general) on the multisummation theorem, that 
any formal isomorphism $\wh F$ between convergent connections may be summed along any non-singular direction to yield a canonical analytic isomorphism.

Let $\pi:\wh \Si\to \IP^1$ denote the real oriented blow-up of $\IP^1$ at $\infty$, 
so $\partial:= \pi^{-1}(\infty)$ is the circle of real  oriented directions emanating from  $\infty\in \IP^1$.
Let $\IH\subset \wh \Si$ be a tubular neighbourhood of the boundary circle $\partial$.
We will call $\IH$ the halo---it is the annulus 
shaded grey on the right of Figure 
\ref{fig: hexagon}.
On $\IH$ we put the rank two local system $\IL_\infty$ of solutions of the formal normal form $A^0$.
This is a graded local system (since $A^0$ is diagonal). 
In the interior of $\Si$ (the white disc in the middle of the figure) we consider the local system
$\cL$ of solutions of the original connection $\nabla$.

Now multisummation yields certain preferred isomorphisms between $\nabla$ and $\nabla^0$, and thus between $\cL$ and $\IL_\infty$. The resulting topological object is the Stokes local system, and can again be axiomatised to classify the original connections.
For this we just need to specify the singular directions and the conditions on the resulting local system.

The singular directions $\IA\subset \partial$ at $\infty$ determined by the irregular type $Q$ are the directions of {\em maximal decay} of one of the ratios $\exp(q_1-q_2)$ or $\exp(q_2-q_1)$ of the exponential factors.
They interlace the Stokes directions and are sometimes called ``anti-Stokes directions''.
This is purely combinatorial:
if $q=\sum_1^k a_iz^i$ then
$\exp(q)$ has maximal decay on the directions where $a_kz^k$ is real and negative.

The multisummation theorem (see e.g. 
\cite{BBRS91})
says that if $d\in \partial\setminus\IA$ is a nonsingular direction then $\wh F$ determines a preferred analytic isomorphism $F_d$ between $A$ and $A^0$ in a sectorial neighbourhood of $d$, 
and these isomorphisms fit together as $d$ varies provided $d$ does not cross a singular direction. 
(In fact in the present set-up only $k$-summation is needed rather than more general multisummation.)

To encode this, let $\wt \Si$ be the surface obtained by deleting a point $e(d)\in \wh \Si$ on the interior boundary $\partial' \neq \partial$ of $\IH$ for each singular direction $d\in \IA$.
(These ``tangential punctures'' $e(d)$ are the small white circles in the figure.)
Thus on each component  
$U\subset \partial'\setminus\{e(d)\st d\in \IA\}$ 
of the interior circle there is a preferred analytic isomorphism $F_U$ 
between $\nabla^0$ and $\nabla$.
Thus we can glue the local systems $\IL_\infty$ and $\cL$ along these components to yield a rank two local system $\IL$ on $\wt \Si$.

To characterise the properties of $\IL$ note that $\IL\bigl\vert_\IH= \IL_\infty$ is graded (by the exponents ${q_1,q_2}$), and for each 
singular direction $d\in \IA$ 
there is a definite dominance ordering of the 
exponential factors: either $\exp(q_1-q_2)$ has maximal decay or $\exp(q_2-q_1)$ does.
This can be encoded by saying ``the root (12) supports $d$'' or  ``the root (21) supports $d$'', respectively (see \cite{bafi} for 
more on the ``rootiness'' of Stokes data).
In turn we can consider the corresponding subgroups
$$
\bmx 1 & * \\ 0 & 1 \emx,
\qquad\text{respectively}\qquad
\bmx 1 & 0 \\ * & 1 \emx 
$$
of automorphisms of the fibre of $\IL$ at 
$d\in \partial$.
These are the Stokes groups $\ISto_d$ for 
$d\in \IA$.
(Intrinsically these groups correspond to adding multiples of the smaller solution on to the larger one, using the grading of $\IL$ over $\IH$.)

Now let $\ga_d$ be a small loop in $\wt \Si$
around $e(d)$ based at $d\in \IA\subset \partial$.
The key fact then is that the monodromy of 
$\IL$ around $\ga_d$ is in $\ISto_d$, and this characterises the local systems on $\wt \Si$ that arise from irregular connections in this way.

\begin{defn}
A Stokes local system for $\Si$ is a rank two local system $\IL$ on $\wt \Si$ such that i) the restriction of $\IL$ to $\IH$ is graded (by $\{q_1,q_2\}$)  and ii) the monodromy of 
$\IL$ around $\ga_d$ is in $\ISto_d$ for all singular directions $d\in \IA$. 
\end{defn}

In turn it is easy to classify Stokes local systems. For example suppose we choose a basepoint $b\in \IH$ and consider the set $\cB$ of isomorphism classes  of Stokes local systems together with a framing of the fibre at $b$ (i.e. a graded isomorphism $\IL_b\cong \IC^2$).
Choosing suitable loops 
generating $\pi_1(\wt \Si,b)$ and taking the 
monodromy then yields the description of $\cB$ given in \eqref{eq: bdescrn}.
Changing the framing corresponds to the action of $T$, and the moment map is the monodromy around $\partial$ (the formal monodromy of $\nabla$).

Finally we can revisit the proof of 
Prop. \ref{prop: sib.ss-sls}.
Given a connection $\nabla$ 
(and thus the data $A,A^0,\wh F$) as above we get both a Stokes filtration and a Stokes local system, and want to relate them. 
Let $\Sect_U$ be a sector at $\infty$ bounded by two consecutive singular directions, and let 
$d$ be the (unique) 
Stokes direction in this sector. 
Now choose a fundamental solution 
$z^\La e^Q$ of $\nabla^0$ on $\Sect_U$ (it is well defined upto right multiplication by a constant diagonal matrix).
Let
$\Phi:= F_Uz^\La e^Q=(v_1,v_2)$ be the fundamental solution of $\nabla$ obtained 
by summing $\wh F$ on $\Sect_U$.
The key relation is now immediate:

\begin{lem}
The column $v_1$ of $\Phi$ is a subdominant solution on one side of $d$ and $v_2$ is subdominant on the other side.
\end{lem}

(More precisely $v_1$ is subdominant on the side of $d$ where $\exp(q_1)$ decays.)
Thus we see how to relate the bases $\Phi$ used to define the canonical Stokes matrices with the bases $\Psi$ used by Sibuya, whose columns are consecutive subdominant solutions (and e.g. explain why we needed to reorder some of the columns in the proof of Prop. \ref{prop: sib.ss-sls}).

\begin{figure}[ht]
	\centering
	\input{slocss2.pstex_t}
	\caption{Stokes local system from Stokes filtered local system}\label{fig: sloc-ss}
\end{figure}

More topologically, if we start with a Stokes local system then near each singular direction we can continue the subdominant graded piece (in the halo) into the interior passing either side of the tangential puncture, to give a well defined ray in the interior, and thus a Stokes filtration.
Conversely given a Stokes filtered local system in the interior we can take the associated graded local system in the halo, and then use the natural projection maps to obtain a Stokes local system. See Figure \ref{fig: sloc-ss}, where all the maps are isomorphisms (but note that not everything is globally defined). 
This is local at the singularity (at infinity here) so we could have an arbitrary curve in the interior (glued to the central circle drawn here).

\renewcommand{\baselinestretch}{1}              %
\normalsize
\bibliographystyle{amsalpha}    \label{biby}
\bibliography{../thesis/syr}

\vspace{0.2cm} 
\noindent
D\'epartement de Math\'ematiques, \\
B\^{a}timent 425, \\
Facult\'e des Sciences d'Orsay \\
Universit\'e Paris-Sud \\
F-91405 Orsay Cedex \\

\noindent
Philip.Boalch@math.u-psud.fr

\end{document}